\documentclass{compositio}
\usepackage{xypic}
\usepackage{eucal}

\usepackage{amsmath}
\usepackage{amstext}
\usepackage{amssymb}
\usepackage{amsthm}
\usepackage{amscd}

\usepackage{lscape}
\usepackage{longtable}

\usepackage{epsfig}
\input txdtools
\let\et=\etexdraw
\def\etexdraw{\drawbb\et}

\theoremstyle{plain}
\newtheorem{thm}{Theorem}[section]
\newtheorem{thm*}{Theorem}

\newtheorem{prop}[thm]{Proposition}
\newtheorem{prop*}[thm*]{Proposition}

\theoremstyle{definition}
\newtheorem{defi}[thm]{Definition}

\theoremstyle{remark}

\DeclareMathOperator{\height}{ht}
\DeclareMathOperator{\Hom}{Hom}

\DeclareMathOperator{\Ext}{Ext}

\DeclareMathOperator{\Ann}{ann}
\DeclareMathOperator{\Nil}{Nil}

\DeclareMathOperator{\grAnn}{gr-ann}
\DeclareMathOperator{\HH}{H}

\newcommand{\twiddle}[2]{I_{#2}\left(#1\right)}

\begin{document}

\def\Ker {\operatorname{Ker}\nolimits}
\def\coker {\operatorname{Coker}\nolimits}
\def\Im {\operatorname{Im}\nolimits}
\def\Image {\operatorname{Image}\nolimits}
\def\Syz {\operatorname{Syz}\nolimits}
\def\initial {\operatorname{in}\nolimits}
\def\Gin{\operatorname{Gin}\nolimits}
\def\Spec{\operatorname{Spec}\nolimits}
\def\D{\operatorname{D}\nolimits}
\def\V{\operatorname{V}\nolimits}
\def\H {\operatorname{H}\nolimits}
\def\E {\operatorname{E}\nolimits}
\def\V {\operatorname{V}\nolimits}
\def\nE {\operatorname{e}\nolimits}
\def\nV {\operatorname{v}\nolimits}
\def\C {\operatorname{\cal C}\nolimits}
\def\row {\operatorname{row}\nolimits}
\def\column {\operatorname{column}\nolimits}

\title
[Parameter test ideals of Cohen Macaulay rings]
{Parameter test ideals of Cohen Macaulay rings}
\author{Mordechai Katzman}
\email{M.Katzman@sheffield.ac.uk}
\address{Department of Pure Mathematics,
University of Sheffield, Hicks Building, Sheffield S3 7RH, United Kingdom\\
{\it Fax number}: +44-(0)114-222-3769}




\begin{abstract}
We describe an algorithm for computing parameter-test-ideals in certain local Cohen-Macaulay rings.
The algorithm is based on the study of a Frobenius map on the injective hull of the residue field of the ring and on the application
of Rodney Sharp's notion of ``special ideals''.

Our techniques also provide an algorithm for computing indices of nilpotency of Frobenius actions on
top local cohomology modules of the ring and on the injective hull of its residue field.
The study of nilpotent elements on injective hulls of residue fields also yields a great simplification of the proof of the
celebrated result in \cite{ABL}.
\end{abstract}

\maketitle

\section{Introduction}\label{Section: Introduction}
This paper deals with various notions originating from the theory of tight closure, which we now review briefly.
Let $S$ be a commutative ring of prime characteristic $p$;
for each ideal $J\subseteq S$ we define the \emph{$e^\text{th}$ Frobenius power} of $J$, denoted $J^{[p^e]}$,
to be the ideal of $S$ generated by $\{ a^{p^e} \,|\, a\in J\}$.
For any ideal $J\subseteq S$ we can then define its \emph{tight closure,} denoted $J^*$,
to be the set of all $a\in S$ such that for some $c\in S$ not in a minimal prime of $S$ we have $c a^{p^e} \in J^{[p^e]}$ for all $e\gg 0$.
Tight closure is indeed a closure operation, in the sense that $J\subseteq J^*$ and $J^{**}=J^*$;
we refer the reader to the seminal paper \cite{Hochster-Huneke-0} and to \cite{Huneke} for a description of tight closure and its properties.

Tight closure has played an important role in many recent advances in Commutative Algebra. A short sample of its useful applications could include
short proofs to some of the
homological conjectures, the study of singularities, positive-characteristic analogues of multiplier ideals and many more.

Among the most interesting and useful results obtained early in the development of the theory of tight closure is the existence of \emph{test-elements}
(cf.~Chapter 2 in \cite{Huneke}). Notice that the element $c\in S$ occurring in the definition of tight closure could depend on
the ideal $J$ and the element $a$. Test-elements are elements $c\in S$ not in any minimal prime such that for \emph{all} ideals $J\subseteq S$ and \emph{all} $a\in S$,
$a\in J^*$ if and only if $c a^{p^e} \in J^{[p^e]}$ for all $e\geq 0$.
Notice, for example, that $J^*=J$ for all ideals $J\subseteq S$ if and only if $1$ is a test-element
(in this case we refer to tight closure as being a \emph{trivial operation}).
Test-elements exist in many rings of interest
(e.g., reduced algebras of finite type over excellent local rings) and they play a vital role. One also defines the
\emph{test-ideal} of $S$ to be the ideal generated by all test-elements.

In many applications one restricts one's attention to local rings and to the tight-closure of ideals generated by systems of parameters.
One then naturally considers the notion of \emph{parameter-test-ideals}: these are
elements $c\in S$ not in any minimal prime such that for all ideals $J\subseteq S$ \emph{generated by a system of parameters} and all $a\in S$,
$a\in J^*$ if and only if $c a^{p^e} \in J^{[p^e]}$ for all $e\geq 0$. It is worth noting that when $S$ is a Gorenstein ring, the notions of
parameter-test-ideals and test-ideals coincide (cf.~Chapter 2 in \cite{Huneke}).

The calculation of tight closure is notoriously hard-- no general algorithm is known and specific calculations are carried out with technical
ad-hoc methods (for example, see \cite{Brenner-Katzman} for such a calculation of a seemingly simple example which settled a major conjecture).
There is not even an algorithm for deciding whether the tight closure operation in a given ring is trivial.

The main aim of this paper is to provide a description of parameter test ideals
of local Cohen-Macaulay rings of prime characteristic $p$. The nature of this description will be
such that it will allow us to give an algorithm for producing these ideals. As a result one also obtains an algorithm
for deciding whether a ring is \emph{$F$-rational,} i.e., whether the tight closure of ideals generated by systems of parameters is trivial;
in the Gorenstein case this property is equivalent to the tight closure of \emph{all} ideals being trivial.

The results in this paper will follow from an analysis of Frobenius
maps on injective hulls of the residue fields $E$ of the ring $S$ under
consideration, i.e., of additive maps $f: E \rightarrow E$ which satisfy $f(s m)=s^p f(m)$ for all $m\in E$ and $s\in S$.
This analysis is inspired by Gennady Lyubeznik's work
on $F$-modules and indeed, a crucial tool used here, namely, the
functors $\Delta^e$ defined in section \ref{Section: A duality}
below, are nothing but ``the first step'' in the construction of
Lyubeznik's $\mathcal{H}$ functors in section 4 of \cite{Lyubeznik}.

The study of $S$-modules with Frobenius maps can be elucidated by treating them as left modules over a certain
skew polynomial ring $S[T; f]$. A crucial ingredient in this paper is Rodney Sharp's recent study of these
modules in general, and of the $S[T; f]$-module structure of the top local cohomology module in particular. In \cite{Sharp}
the parameter test ideal of $S$ was described in terms of certain $S[T; f]$-submodules of the top local cohomology of $S$, and it is this description on which our explicit description and algorithm is based on.

Along the way we gain new insights into the $S[T; f]$-module structure of injective hulls of residue fields which translate
into new results. One such result is an algorithm for computing the index of nilpotency (in the sense of
section 4 of \cite{Lyubeznik}) of top local cohomology modules, which, together with the results in
\cite{Katzman-Sharp} translate into an algorithm for computing the Frobenius closure of
parameter ideals in Cohen-Macaulay local rings and in view of  \cite{Huneke-Katzman-Sharp-Yao}
provide an important ingredient for the corresponding computation in generalized Cohen-Macaulay rings as well.
Another spinoff is a very simple proof of a crucial ingredient in
\cite{ABL} which together with Corollary 3.6 there gives an alternative proof of the fact that
for a power series ring $R$ of prime characteristic, for all nonzero $f\in R$, $1/f$
generates $R_f$ as a $D_R$-module.

\section{Frobenius maps}\label{Section: Frobenius maps}

Let $S$ be a commutative ring of prime characteristic $p$ and let $M$ be an $S$-module.
A \emph{Frobenius map} on $M$ is a $\mathbb{Z}$-linear map $\phi: M \rightarrow M$ with the property that
$\phi(s m)=s^p \phi(m)$ for all $s\in S$ and $m\in M$.
The fundamental example of a Frobenius map is \emph{the} Frobenius map  $f: S \rightarrow S$
given by $f(s)= s^p$. The Frobenius map allows us to endow $S$ with a structure of an $S$-bimodule:
as a left $S$-module it has the usual $S$-module structure whereas $S$ acts on itself on the right
via the Frobenius map. We shall denote this bimodule $F_S(S)$ and so for all $a\in F(S)$ and $s\in S$,
$s \cdot a=sa$ while $a\cdot s=s^p a$, where $\cdot$ denotes the action of $S$.
We can extend this construction to obtain the \emph{Frobenius functor} $F$
sending any $S$-module $M$ to $F_S(M)=F_S(S)\otimes_S M$ where $S$ acts on $F_S(M)$ via its left-action on
$F_S(S)$, so for $s\otimes m\in F_S(M)$ and $t\in S$
we have $t\cdot (s\otimes m)=ts\otimes m$ and $(s\otimes tm)=s \cdot t \otimes m=t^p s \otimes m$.
We shall repeatedly (and tacitly) use the fact that the functor $F_S$ is exact whenever $S$ is regular (Theorem 2.1 in \cite{Kunz}).

Iterations  $\phi^e=\underbrace{\phi \circ\dots \circ \phi }_{e \text{ times}}$
of Frobenius maps  $\phi: M \rightarrow M$ result in maps
$\phi^e: M \rightarrow M$ which satisfy
$\phi(s m)=s^{p^e} \phi(m)$ for all $s\in S$ and $m\in M$.
More generally, we will consider the set $\mathcal{F}^e(M)$ of all
$\mathbb{Z}$-linear maps $\psi: M \rightarrow M$ which satisfy
$\psi(s m)=s^{p^e} \psi(m)$ for all $s\in S$ and $m\in M$.
We can give $\mathcal{F}^e(M)$
the structure of an $S$-module: for any $\psi\in \mathcal{F}^e(M)$ and  $a\in S$
we simply let $a \psi$ the map sending $m\in M$ to $a \psi(m)$.
Furthermore, we can define a product in $\mathcal{F}(M):= \mathcal{F}^e(M)$ to be composition of
$\mathbb{Z}$-linear maps and thus endow $\mathcal{F}(M)$ with the structure of an $S$-algebra.

The iteration of the Frobenius map on $R$ leads one
to the iterated Frobenius functors $F^i_R(-)$ which are defined
for all $i\geq 1$ recursively by $F^1_R(-)=F_R(-)$ and
$F^{i+1}_R(-)=F_R\circ \left(F^i_R(-)\right)$ for all $i\geq 1$.
These higher Frobenius functors are also exact whenever $S$ is regular.

In this paper we will be interested in studying Frobenius maps on injective hulls of residue fields
and top local cohomology modules. An example of the latter
when $S$ is local and $d$-dimensional can be obtained as follows.
The top local cohomology module $\HH^d_\mathfrak{m}(S)$ can be computed as the direct limit of
$$
\frac{S}{(x_1, \dots, x_d)S} \xrightarrow[]{x_1 \cdot \ldots \cdot x_d} \frac{S}{(x_1^2, \dots, x_d^2)S}  \xrightarrow[]{x_1 \cdot \ldots \cdot x_d} \dots
$$
where $x_1, \dots, x_d$ is a system of parameters of $S$.
We can define a Frobenius map $\phi\in\mathcal{F}^e\left( \HH^d_\mathfrak{m}(S)\right)$ on this direct limit by mapping the coset
$a + (x_1^n, \dots, x_d^n)S$ in the $n$-th component of the direct limit
to the coset $a^{p^e} + (x_1^{np^e}, \dots, x_d^{np^e})S$
in the $np^e$-th component of the direct limit.
It is not hard to verify that this is indeed a well defined map from $\HH^d_\mathfrak{m}(S)$ to $\HH^d_\mathfrak{m}(S)$
and that it is a Frobenius map.
An important observation used in this paper is
the fact that when $S$ is Cohen-Macaulay,
the Frobenius map $\phi\in \mathcal{F}^1(\HH^d_\mathfrak{m}(S))$ described above generates
the $S$-algebra $\mathcal{F}(\HH^d_\mathfrak{m}(S))$ (cf. Example 3.7 in \cite{Lyubeznik-Smith}).

\bigskip
A different and fruitful way of thinking about Frobenius maps on $M$ and their iterations is
as left-module structures over certain skew-commutative rings. Given any commutative ring $S$
we can construct a skew commutative ring $S[T; f^e]$ as follows:
as an $S$-module it will be the free module $\displaystyle \oplus_{i=0}^\infty S T^{i}$
and we extend the rule $Ts=s^{p^e}T$ for all $s\in S$ to a (non-commutative!) multiplicative structure on $S[T; f^e]$.
Given a Frobenius map $\phi\in \mathcal{F}^e(M)$ on an $S$-module $M$, we can then turn it into a left
$S[T; f^e]$-module by extending the rule $T m=\phi(m)$ for all $m\in M$.
The fact that this gives $M$ the structure of a left $S[T; f^e]$-module is simply because
for all $s\in S$ and $m\in M$,
$$T (s m) = \phi(sm)=s^{p^e} \phi(m)= s^{p^e} T m = (T s) m .$$
This approach has been taken in many previous papers, the most relevant to us being \cite{Sharp}.

\section{A duality}\label{Section: A duality}

In this section we set up the main tool, based on Matlis Duality,
which will enable us to explore
$R[T; f]$-module structures of certain Artinian modules.

\emph{
Henceforth in this paper $(R, \mathfrak{m})$ will denote a complete local regular ring of characteristic $p$.
We shall denote the injective hull of $R/\mathfrak{m}$ with $E$ and $(-)^\vee$ shall denote the functor
$\Hom_R(- , E)$.}

Let $M$ be any $R$-module and for all $m\in M$, let $e_m\in M^{\vee\vee}$ be defined by
$e_m(g)=g(m)$ for all $g\in M^\vee$.
Matlis duality states that for all $R$-modules $M$ which are either Noetherian or Artinian,
the map $M  \rightarrow M^{\vee\vee}$ which sends $m\in M$ to $e_m$ is an isomorphism of $R$-modules.
If now $M$ is an $R[T;f^e]$-module, this map endows $M^{\vee\vee}$ with a structure of an $R[T;f^e]$-module
defined by $T e_m= e_{Tm}$ for all $m\in M$, so now we may identify $M$ and $M^{\vee\vee}$ as $R[T;f^e]$-modules.

Since $R$ is complete, a straightforward modification of Lemma 4.1 in \cite{Lyubeznik}
provides us with a natural, functorial isomorphism $\gamma^e_M: F_R^e(M)^\vee \rightarrow F^e_R(M^\vee)$
defined for all Artinian $R$-modules. We shall use this isomorphism repeatedly in this section.

Fix now an ideal $I\subseteq R$ and write $S=R/I$.
Let $\mathcal{C}^e$ be the category of Artinian $S[T;f^e]$-modules.
Let
$\mathcal{D}^e$ be the category of $R$-linear maps $M \rightarrow F^e_R(M)$ where $M$ is
a finitely generated $S$-module and where a morphism between $M\xrightarrow[]{a} F_R^e(M)$ and
$N\xrightarrow[]{b} F_R^e(N)$ is a commutative diagram of $S$-linear maps
\begin{equation*}
\xymatrix{
M \ar@{>}[d]^{a} \ar@{>}[r]^{\mu} & N \ar@{>}[d]^{b}\\
F_R^e(M) \ar@{>}[r]^{F^e_R(\mu)} & F_R^e(N)\\
}
\end{equation*}

In this section we construct a pair of functors
$\Delta^e: \mathcal{C}^e \rightarrow \mathcal{D}^e$ and
$\Psi^e: \mathcal{D}^e \rightarrow \mathcal{C}^e$
in a way that for all $M\in \mathcal{C}^e$,
the $S[T; f^e]$-module $\Psi^e \circ \Delta^e (M)$ is canonically isomorphic to $M$
and for all $D=(N\xrightarrow[]{u} F^e_R(N))\in \mathcal{D}^e$,
$\Delta^e \circ \Psi^e (D)$ is canonically isomorphic to $D$.

The functor $\Delta^e$ just the ``first step'' in the construction of Gennady Lyubeznik's functor
$\mathcal{H}_{R,S}$, i.e.,
for $M\in \mathcal{C}^e$ we have an $R$-linear map
$\alpha_M: F^e_R(M) \rightarrow M$
given by $\alpha(r\otimes m)=r T m$ for all $r\in R$ and $m\in M$.
Applying $(-)^\vee$  to the map $\alpha$ one obtains an $R$-linear map
$\alpha_M^\vee: M^\vee \rightarrow  F_R^e(M)^\vee$.
We now define the $\Delta(M)$ to be the map
$M^\vee \xrightarrow[]{\gamma_M \circ \alpha_M^\vee} F_R^e(M^\vee)$.

To define  $\Psi^e$ we retrace the steps above;
given a finitely generated $S$-module $N$ and a $R$-linear map
$a: N \rightarrow F^e_R(N)$
we define $\Psi^e(-)$ to coincide with the functor $(-)^\vee$ as a functor of $S$ modules
giving $\Psi^e(N)$ the additional structure of an $S[T;f^e]$ module structure as follows.

We apply ${}^\vee$ to the map $a$ above to obtain a map
$a^\vee : F_R^e(N)^\vee \rightarrow N^\vee$.
We next obtain a map $\epsilon: F_R^e\left(N^\vee\right) \rightarrow F_R^e(N)^\vee$ as the following composition:
$$F_R^e\left(N^\vee\right) \cong
F_R^e\left(N^\vee\right)^{\vee\vee} \xrightarrow[]{\left(\gamma_{N^\vee}^\vee\right)^{-1}}
F_R^e\left(N^{\vee\vee}\right)^{\vee} \cong
F_R^e\left(N\right)^{\vee}.$$
We now obtain a functorial map
$b=a^\vee\circ \epsilon: F^e_R(N^\vee) \rightarrow N^\vee$ and
we define the action of $T$ on $N^\vee$ by defining
$T n=b(1\otimes n)$  for all $n\in N^\vee$.

\begin{thm}\label{Theorem: the two functors}
The functors
$\Delta^e: \mathcal{C}^e \rightarrow \mathcal{D}^e$ and
$\Psi^e: \mathcal{D}^e \rightarrow \mathcal{C}^e$ are exact.
For all $M\in \mathcal{C}^e$,
the $S[T; f]$-module $\Psi^e \circ \Delta^e (M)$ is canonically isomorphic to $M$.
For all $D=(N\xrightarrow[]{u} F^e_R(N))\in \mathcal{D}^e$,
$\Delta^e \circ \Psi^e (D)$ is canonically isomorphic to $D$.
\end{thm}

\begin{proof}
The exactness of the functors follows from the exactness of the functors $\Hom_R(-, E)$ and $F^e_R$.

To prove the second statement we notice that for all $M\in \mathcal{C}^e$,
$\Psi^e \circ \Delta^e (M)$ is $M^{\vee\vee}$ which we identify as an $S$-module with $M$ by identifying each $m\in M$
with the $e_m\in M^{\vee\vee}$ which we defined at the beginning of this section.
We want to show that
this identification is an isomorphism of $S[T;f^e]$-modules, and to do so we now describe
$T e_m$ for all $e_m\in\Psi^e \circ \Delta^e (M)$. This will be the image of $1\otimes e_m$ under the map
\begin{equation}\label{eqn1}
F_R^e\left(M^{\vee\vee}\right) \xrightarrow[]{i_1}
F_R^e\left(M^{\vee\vee}\right)^{\vee\vee} \xrightarrow[]{\left(\gamma_{M^{\vee\vee}}^\vee\right)^{-1}}
F_R^e\left(M^{\vee\vee\vee}\right)^{\vee} \xrightarrow[]{i_2}
F_R^e\left(M^{\vee}\right)^{\vee} \xrightarrow[]{\gamma_M^\vee}
F_R^e(M)^{\vee\vee} \xrightarrow[]{\alpha^{\vee\vee}}
 M^{\vee\vee}
\end{equation}
where $i_1, i_2$ are the isomorphisms induced from the isomorphism of functors $(-)\cong(-)^{\vee\vee}$.
The functoriality of $\gamma_{(-)}$ implies that we have a commutative diagram
\begin{equation*}
\xymatrix{
F_R^e\left(M^{\vee\vee\vee}\right)^{\vee} \ar@{>}[d]^{\gamma_{M^{\vee\vee}}^\vee} \ar@{>}[r]^{i_2} &
F_R^e\left(M^{\vee}\right)^{\vee}  \ar@{>}[d]^{\gamma_M^\vee}\\
F_R^e(M^{\vee\vee})^{\vee\vee} \ar@{>}[r]^{i_1^{-1}} &
F_R^e(M)^{\vee\vee} \\
}
\end{equation*}
We may now rewrite the composition in (\ref{eqn1}) as
\begin{equation*}
F_R^e\left(M^{\vee\vee}\right) \xrightarrow[]{i_1}
F_R^e\left(M^{\vee\vee}\right)^{\vee\vee} \xrightarrow[]{\left(\gamma_{M^{\vee\vee}}^\vee\right)^{-1}}
F_R^e\left(M^{\vee\vee\vee}\right)^{\vee} \xrightarrow[]{\gamma_{M^{\vee\vee}}^\vee}
F_R^e(M^{\vee\vee})^{\vee\vee}  \xrightarrow[]{i_1^{-1}}
F_R^e(M)^{\vee\vee} \xrightarrow[]{\alpha^{\vee\vee}}
 M^{\vee\vee}
\end{equation*}
which simplifies into
$$F_R^e(M)^{\vee\vee} \xrightarrow[]{\alpha^{\vee\vee}}  M^{\vee\vee}.$$

If we now start with $D=(N \xrightarrow[]{u} F^e_R(N))\in \mathcal{D}^e$
$$\Delta^e\circ \Psi^e(D): N^{\vee\vee} \rightarrow F^e \left(N^{\vee\vee}\right)$$
is given by the composition
\begin{equation}\label{eqn2}
N^{\vee\vee} \xrightarrow[]{a^{\vee\vee}}
F_R^e\left(N\right)^{\vee\vee} \xrightarrow[]{i_3}
F_R^e\left(N^{\vee\vee}\right)^{\vee\vee} \xrightarrow[]{\left(\left(\gamma_{N^\vee}^\vee\right)^{-1}\right)^\vee}
F_R^e\left(N^{\vee}\right)^{\vee\vee\vee} \xrightarrow[]{i_4}
F_R^e\left(N^{\vee}\right)^\vee \xrightarrow[]{\gamma_{N^\vee}}
F_R^e\left( N^{\vee\vee} \right)
\end{equation}
where $i_1, i_2$ are the isomorphisms induced from the isomorphism of functors $(-)\cong(-)^{\vee\vee}$.
Now $\left(\left(\gamma_{N^\vee}^\vee\right)^{-1}\right)^\vee=\left(\gamma_{N^\vee}^{\vee\vee}\right)^{-1}$
and the functoriality of $\gamma_{(-)}$ implies that we have a commutative diagram
\begin{equation*}
\xymatrix{
F_R^e\left(N^{\vee\vee}\right)^{\vee\vee}  \ar@{>}[d]^{i_3^{-1}} \ar@{>}[r]^{\left(\gamma_{N^\vee}^{\vee\vee}\right)^{-1}} &
F_R^e\left(N^{\vee}\right)^{\vee\vee\vee}  \ar@{>}[d]^{i_4}\\
F_R^e(N^{\vee\vee}) \ar@{>}[r]^{{\gamma_{N^\vee}}^{-1}} &
F_R^e\left(N^{\vee}\right)^\vee \\
}
\end{equation*}
We may now rewrite the composition in (\ref{eqn2}) as
\begin{equation*}
N^{\vee\vee} \xrightarrow[]{a^{\vee\vee}}
F_R^e\left(N\right)^{\vee\vee} \xrightarrow[]{i_3}
F_R^e\left(N^{\vee\vee}\right)^{\vee\vee} \xrightarrow[]{i_3^{-1}}
F_R^e(N^{\vee\vee}) \xrightarrow[]{\gamma_{N^\vee}^{-1}}
F_R^e\left(N^{\vee}\right)^\vee \xrightarrow[]{\gamma_{N^\vee}}
F_R^e\left( N^{\vee\vee} \right)
\end{equation*}
which simplifies to
$a^{\vee\vee}$.
\end{proof}

Throughout this paper, when $e=1$ we will drop the subscript $e$ from our notation.
Thus $\mathcal{C}^1=\mathcal{C}$, $\Delta^1=\Delta$, etc.

As mentioned before, the functor $\Delta$ is a building block for another functor described in section 4 of \cite{Lyubeznik}.
This functor, denoted with $\mathcal{H}_{R,S}$, is a functor from $\mathcal{C}$ to the category of $F$-finite $F_R$-modules
(see sections 1,2 and 3 in  \cite{Lyubeznik} for definition and properties)
and is obtained as follows.
For $M\in \mathcal{C}$ write $\Delta(M)=\left(N\xrightarrow[]{u} F_R(N)\right)$.
Now $\mathcal{H}_{R,S}(M)$ is defined to be the direct limit of
$$ N\xrightarrow[]{u} F_R(N) \xrightarrow[]{F_R(u)} F^2_R(N) \xrightarrow[]{F^2_R(u)} \dots .$$
Various useful properties of Lyubeznik's functor can be found in section 4 of \cite{Lyubeznik}.

\section{Frobenius maps on injective hulls}\label{Section: Frobenius maps on injective hulls}

\emph{Henceforth in this paper we shall fix an ideal $I\subseteq R$ and denote $R/I$ with $S$.}

In this section we will first apply the tools developed in section \ref{Section: A duality}
to yield a  description of possible $S[T; f^e]$-module structures of $E_S=E_S(S/\mathfrak{m}S)$, the injective hull of
the residue field of $S$.
This description is not new: it is contained in Proposition 5.2 in \cite{Lyubeznik-Smith}.
Later in the section we shall use this description to describe explicitly the nilpotent elements of $E_S(S/\mathfrak{m}S)$.

\begin{prop}
The $S[T; f^e]$-module structures on $E_S$ are given by
$$\Psi^e\left( R/I \xrightarrow[]{g}  R/I^{[p^e]}\right)$$
where the map above is given by multiplication by some $g\in \left(I^{[p^e]} : I\right)$.
\end{prop}
\begin{proof}
Clearly, all $R$ linear maps $ R/I\rightarrow  R/I^{[p^e]}$ are given by multiplication by some $g\in \left(I^{[p^e]} : I\right)$.
The proposition now follows from Theorem \ref{Theorem: the two functors} and the fact that $E^{\vee}\cong S$.
\end{proof}

\bigskip
The bijection between $R$-linear maps $R/I \xrightarrow[]{}  R/I^{[p^e]}$ and
$S[T; f^e]$-module structures on $E_S$ has been described explicitly in Chapter 3 of \cite{Blickle} as follows.
First, notice that $E$, thought of as the direct limit of
$$\frac{R}{(y_1, \dots, y_n)R} \xrightarrow[]{y_1\cdot\ldots\cdot y_n}
\frac{R}{(y_1^2, \dots, y_n^2)R} \xrightarrow[]{y_1\cdot\ldots\cdot y_n} \dots$$
where $y_1, \dots, y_n$ is a system of parameters for $R$,
has a natural Frobenius map given by
$$\phi\left( r + (y_1^s, \dots, y_n^s) \right)= r^p + (y_1^{s p}, \dots, y_n^{s p}) \in \frac{R}{(y_1^{s p}, \dots, y_n^{s p})R} .$$
Now if  $u\in (I^{[p]}:_R I)$, $u \phi$, which is also a Frobenius map on $E$, will restrict
to a Frobenius map on $E_S=\Ann_E I$ because for all $m\in \Ann_E I$,
$$I u \phi(m) \subseteq I^{[p]} \phi(m)= \phi(I m)= \phi(0)=0 .$$
In Chapter 3 of \cite{Blickle} it is shown that \emph{all} Frobenius maps on $E_S$ are obtained in this way.

\bigskip
Henceforth in this section we shall assume that $E_S$ has a given $S[T; f]$-module structure.
Our next aim is to describe the $S[T; f]$-submodules of $E_S$.
Later in the section we shall use this description to describe explicitly the nilpotent elements of $E_S$.
We start by recalling that the set of $S$-submodules of $E_S$ is
$\displaystyle \{ \Ann_{E_S} J \,|\, J\subseteq S\}$ (cf. Theorem 5.21 in \cite{Sharpe-Vamos}).
If we now asked for a description of the $S[T;f]$-submodules of $E_S$,
the answer would obviously be ``all $\Ann_{E_S} J$ which happen to be $S[T;f]$-submodules of $E_S$''.
With this in mind we define the following.

\begin{defi}
An ideal $J\subseteq S$ is called an $E_S$-ideal if
$\displaystyle \Ann_{E_S} J$
is an $S[T;f]$-submodule of $E_S$.
An ideal $J\subseteq R$ is called an $E_S$-ideal if it contains $I$ and its image in $S$ is an $E_S$-ideal.
\end{defi}
Notice that for an ideal $J\subseteq S$, being an $E_S$-ideal is equivalent to
$\displaystyle \Ann_{E_S} J= \Ann_{E_S} J S[T;f]$. We also note that when $S$ is Gorenstein the notion
of $E_S$ ideals coincides with that of $F$-ideals studied in \cite{Smith2}.

\begin{thm}
Let $u\in R$ be such that
$\Delta(E_S)$ is the map $\displaystyle \frac{R}{I} \xrightarrow[]{u} \frac{R}{I^{[p]}}$.
The $E_S$ ideals in $R$ consist of all ideals $ L \subseteq R$ containing $I$ for which
$u L \subseteq L^{[p]}$.
\end{thm}
\begin{proof}
Assume first that $L$ is an $E$-ideal.
Apply the functor $\Delta$
to the short exact sequence of $S[T;f]$ modules
\begin{equation}\label{eqn2a}
0 \rightarrow \Ann_{E_S} L \rightarrow E_S \rightarrow E/\Ann_{E_S} L \rightarrow 0
\end{equation}
to obtain the following short exact sequence in $\mathcal{D}$
\begin{equation}\label{CD1}
\xymatrix{
0 \ar@{>}[r]^{} &
\displaystyle\frac{L}{I} \ar@{>}[r]^{} \ar@{>}[d]^{u} &
\displaystyle\frac{R}{I} \ar@{->}[r]^{} \ar@{>}[d]^{u} &
\displaystyle\frac{R}{L} \ar@{>}[r]^{} \ar@{>}[d]^{u}  &
0 \\
0 \ar@{>}[r]^{} &
\displaystyle \displaystyle\frac{L^{[p]}}{I^{[p]}}  \ar@{>}[r]^{} &
\displaystyle \displaystyle\frac{R}{I^{[p]}}        \ar@{>}[r]^{} &
\displaystyle \displaystyle\frac{R}{L^{[p]}}        \ar@{>}[r]^{} &
0 \\
}
\end{equation}
and we must have
$u L \subseteq L^{[p]}$.

On the other hand, if $u L \subseteq L^{[p]}$,
we can construct the commutative diagram (\ref{CD1}), and an application
of the functor $\Psi$ gives back the short exact sequence (\ref{eqn2a})
and we deduce that $L$ is an $E_S$-ideal.
\end{proof}

\bigskip
We now
turn our attention to the nilpotent elements of $E_S$, i.e., to the $S[T; f]$-submodule of $E_S$
$$\Nil(E_S)=\left\{ m\in E_S \,|\, T^e m=0 \text{ for some } e\geq 0 \right\} . $$
Recall that we can write   $\Nil(E_S)$ as  $\Ann_{E_S} J S[T;f]$ for some $E_S$-ideal $J\subseteq R$.
Also, it is known that there exists an $\eta\geq 1$ such that $T^\eta \Nil(E_S)=0$
(cf.~\cite[Proposition 1.11]{Hartshorne-Speiser} and \cite[Proposition 4.4]{Lyubeznik}).
This invariant of $S$ plays an important role in the study of the Frobenius closure
(see \cite{Katzman-Sharp} and \cite{Huneke-Katzman-Sharp-Yao}).
We now describe the ideal $J$ and the index of nilpotency $\eta$.

\begin{defi}
For all $e\geq 1$ write $\nu_e=1+p+\dots+p^{e-1}$.
\end{defi}

\begin{prop}
Let the map $\Delta(E_S)=(R/I \rightarrow R/I^{[p]})$ be given by multiplication by $u\in R$.
Consider $E_S$ as an $S[\Theta; f^e]$-module where for all $m\in E_S$, we define $\Theta m=T^e m$.
The map $\Delta^e(E_S)=(R/I \rightarrow R/I^{[p^e]})$ is given by multiplication by $u^{\nu_e}$.
\end{prop}
\begin{proof}
For all $e\geq 1$ the $R$-linear map $\alpha: F^{e}(E_S) \rightarrow E_S$ defined by
$\alpha( r \otimes m )= r T^{e} m$ can be factored as
$\alpha=\alpha_1 \circ \dots \circ \alpha_e $
where for all $1\leq i\leq e$, $\alpha_i$ is the $R$-linear map $\alpha_i : F^{i}(M) \rightarrow F^{i-1}(M)$
defined by $\alpha_i(r\otimes m)=r\otimes T m$.
Also, it is not hard to see that  $\alpha_{i+1}=F(\alpha_i)$ for all $1\leq i\leq e$.

Now for all $e\geq 1$, the map  $\Delta^e(E_S)=(R/I \rightarrow R/I^{[p^e]})$ is given
by $\gamma^e_{E_S} \circ \alpha^\vee$. It follows from the construction of $\gamma^e_{E_S}$ that,
if we identify $E_S^\vee$ with $R/I$, $\gamma^e_{E_S} : R/I^{[p^e]} \rightarrow R/I^{[p^e]}$ is the identity map.
Now
$$\alpha^\vee=
\alpha_e^\vee \circ \dots \circ \alpha_1^\vee =
F^{e-1} (u) \circ \dots \circ u=
u^{p^{e-1}} \circ \dots \circ u=u^{\nu_e}. $$
\end{proof}

\begin{thm}\label{Theorem: description of nilpotent elements}
Let the map $\Delta(E_S)=(R/I \rightarrow R/I^{[p]})$ be given by multiplication by $u\in R$.
For all $e\geq 1$ let $J_e$ be the smallest ideal of $R$ which contains $I$ and  such that $u^{\nu_e}\in J_e^{[p^e]}$.
There exists an $\alpha\geq 1$ such that $J_\alpha=J_{\alpha+1}$
and for this $\alpha$,
$\Ann_{E_S} J_\alpha$  coincides with the $S[T;f]$-module $\Nil(E_S)$
of nilpotent elements of $E_S$.
Furthermore, the index of nilpotency of $\Nil(E_S)$, if not zero, is the smallest such $\alpha$.
\end{thm}
\begin{proof}
For all $e\geq 1$ let $N_e=\{m\in E_S \,|\, T^e m=0 \}$
and write $N_e=\Ann_{E_S} L_e$ for some
$E_S$-ideal $L_e$.

Notice that $\Delta^e(N_e)=(R/L_e \rightarrow R/L_e^{[p^e]})$ and that the previous proposition
implies that this map is given by multiplication by $u^{\nu_e}$.
It follows from the construction of $\Delta^e(N_e)=(R/L_e \xrightarrow[]{u^{\nu_e}} R/L_e^{[p^e]})$
that this map is the zero map, i.e., $u^{\nu_e}\in L_e^{[p^e]}$; now the minimality of $J_e$
implies that $J_e\subseteq L_e$ and $\Ann_{E_S} L_e \subseteq \Ann_{E_S} J_e$.

On the other hand, the map $R/J_e \xrightarrow[]{u^{\nu_e}} R/J_e^{[p^e]}$ is the zero map
and so $T^e$ kills
$$\Psi^e\left(R/J_e \xrightarrow[]{u^{\nu_e}} R/J_e^{[p^e]}\right)\cong \Ann_{E_S} J_e $$
hence $\Ann_{E_S} J_e \subseteq N_e = \Ann_{E_S} L_e$ and we deduce that
$\Ann_{E_S} J_e = N_e$.

Proposition 4.4 in \cite{Lyubeznik} now implies that  $\Nil(E_S)=N_\alpha$ for some $\alpha\geq 1$,
and since $\Nil(E_S)$ is the union of the ascending chain $\{ N_\alpha \}_{\alpha\geq 1}$, we see that
$N_\beta=N_\alpha$ for all $\beta\geq \alpha$.
Also, if $N_\alpha=N_{\alpha+1}$ but $\Nil(E_S)\neq N_\alpha$, pick any non-zero
$m\in \Nil(E_S)\setminus N_\alpha$ and let $i\geq 0$ be minimal such that $T^i m\notin  N_\alpha$.
Now $T^{\alpha+1} (T^i m)= T^\alpha T^{i+1} m = 0$ so
$T^{i} m\in N_{\alpha+1} \setminus N_\alpha$, a contradiction.
\end{proof}

We shall see in section \ref{Section: The star-closure} how to compute this smallest ideal $J\supseteq I$ for which $u^{\nu_e}\in J^{[p^e]}$.

\bigskip
We conclude this section by exhibiting another ``naturally occurring" $S[T;f]$-submodule of $E_S$.
\begin{thm}\label{Theorem: largest nilpotent quotient}
Let the map $\Delta(E_S)=(R/I \rightarrow R/I^{[p]})$ be given by multiplication by $u\in R$.
For all $\alpha\geq 0$ write
$L_\alpha=\left(I^{[p^\alpha]} :_R u^{\nu_\alpha}\right)$.
\begin{itemize}
    \item [(a)] The sequence of ideals $\{ L_\alpha \}_{\alpha\geq 1}$ is an ascending sequence.
    \item [(b)] If $L_A=L_{A+1}$ then $L_\alpha=L_A$ for all $\alpha\geq A$ and $L_A$ is an $E_S$-ideal.
    \item [(c)] Let $L$ be the stable value of $\{L_\alpha\}_{\alpha\geq 1}$.
    The quotient $E_S/\Ann_{E_S} L$ is nilpotent and for any $E_S$-ideal $K$,
    $E_S/\Ann_{E_S} K$ is nilpotent if and only if $K\subseteq L$.
\end{itemize}
\end{thm}

\begin{proof}
For all $\alpha\geq 1$ the map $g_\alpha: R/I \rightarrow R/I^{[p^\alpha]}$ given by the composition
$$R/I \xrightarrow[]{u} R/I^{[p]} \xrightarrow[]{u^p} R/I^{[p^2]}
\xrightarrow[]{u^{p^2}} \dots  \xrightarrow[]{u^{p^{\alpha-1}}} R/I^{[p^\alpha]} $$
is just the map $g_\alpha : R/I \xrightarrow[]{u^{\nu_\alpha}} R/I^{[p^\alpha]}$
given by multiplication by $u^{\nu_\alpha}$ and whose kernel is $L_\alpha$.
These kernels form an ascending chain, so (a) follows.

The first statement in (b) follows from Proposition 2.3(b) in \cite{Lyubeznik}.
To prove the second statement, we first notice that since
$u I \subseteq I^{[p]}$,
$u^{p^\alpha} I^{[p^\alpha]} \subseteq I^{[p^{\alpha+1}]}$ for all $\alpha\geq 1$
hence $u^{\nu_\alpha} I \subset I^{[p^\alpha]}$ and we deduce that $I\subset L_\alpha$ for all $\alpha\geq 1$.
To show that  $L_A$ is an $E_S$-ideal it remains to prove that $u L_A \subseteq L_A^{[p]}$:
\begin{eqnarray*}
u L_A & = & u L_{A+1} \\
 & \subseteq & \left( I^{[p^{A+1}]} :_R u^{p+p^2+\dots +p^{A}} \right)\\
 & = & \left( \left(I^{[p^{A}]}\right)^{[p]}  :_R \left(u^{N_A}\right)^p \right)\\
 & = & \left( I^{[p^{A}]}  :_R u^{N_A} \right)^{[p]}\\
 & = & L_A^{[p]}\\
\end{eqnarray*}
where the penultimate equality is a consequence of the exactness of $F_R(-)$.

Let $K$ be an $E_S$-ideal for which $E_S/\Ann_{E_S} K$ is nilpotent and choose some
$e\geq 1$ for which $T^e \left(E_S/\Ann_{E_S} K\right)=0$.
An application of $\Delta^e$ to the short exact sequence
$$ 0 \rightarrow \Ann_{E_S} K \rightarrow E_S \rightarrow E_S/ \Ann_{E_S} K \rightarrow 0$$
produces the following short exact sequence in $\mathcal{D}^e$
\begin{equation*}
\xymatrix{
0  \ar@{>}[r]^{} & K/I \ar@{>}[d]^{u^{\nu_e}} \ar@{>}[r]^{} & R/I         \ar@{>}[r]^{} \ar@{>}[d]^{u^{\nu_e}} & R/K \ar@{>}[d]^{u^{\nu_e}} \ar@{>}[r]^{} & 0\\
0  \ar@{>}[r]^{} & K^{[p^e]}/I^{[p^e]}       \ar@{>}[r]^{} & R/I^{[p^e]}  \ar@{>}[r]^{}                        & R/K^{[p^e]}               \ar@{>}[r]^{} & 0\\
}
\end{equation*}
where the leftmost vertical map is the zero map, i.e., $u^{\nu_e} K \subseteq I^{[p^e]}$ and hence
$K\subseteq L_e$.

\end{proof}

\section{The $\star$-closure}\label{Section: The star-closure}

The
statements of Theorems \ref{Theorem: description of nilpotent elements} and \ref{Theorem: largest nilpotent quotient}
in the previous section referred to certain smallest ideals $J\subseteq R$
with the property that $J^{[p^e]}$ contains a given ideal.
The aim of this section is to establish the existence of these ideals and to
describe an algorithm for computing them.

Throughout this section $T$ will denote a Noetherian regular ring of prime characteristic $p$.

\begin{defi}
Let $e\geq 1$.
For any ideal $A\subseteq T$ we define
$$\mathcal{G}^e(A)=\left\{ L \,|\, L \subseteq T \text{ an ideal, } A\subseteq L^{[p^e]} \right\}$$
and
$\displaystyle \twiddle{A}{e}=\bigcap_{L\in \mathcal{G}^e(A)} L$.
\end{defi}

Note that  in general there is no reason why $\twiddle{A}{e}$ should be in $\mathcal{G}^e(A)$.
Recall that a $T$-module $M$ is $\cap$-flat if it is flat and if for all sets of $T$-submodules
$\left\{ N_\lambda \right\}_{\lambda\in\Lambda}$ of a finitely generated module $N$,
$$M\otimes_T \cap_{\lambda\in\Lambda} N_\lambda= \cap_{\lambda\in\Lambda} \left(M\otimes_T N_\lambda\right)$$
(cf.~\cite{Hochster-Huneke-1} pp.~41). Notice that free modules are $\cap$-flat.

\begin{prop}\label{Proposition: properties of widetilde}
Let $e\geq 1$ and assume that $T^{1/p^e}$ is a $\cap$-flat $T$-module.
Let $A\subseteq T$ be an ideal.
\begin{itemize}
    \item[(a)] $\twiddle{A}{e}\in \mathcal{G}^e(A)$ and is the minimal element of $\mathcal{G}^e(A)$.
    \item[(b)] Let $B\subseteq T$ be any ideal.
    The smallest ideal $J\subseteq T$ which contains both $A^{[p^e]}$ and $B$ is $\twiddle{A}{e}+B$.
    \item[(c)] If $A=A_1+\dots+A_s$ then $\twiddle{A}{e}=\twiddle{A_1}{e}+\dots+\twiddle{A_s}{e}$.
\end{itemize}
\end{prop}

\begin{proof}
The first statement is an immediate consequence of the fact that the $T$-module $T^{1/{p^e}}$
is assumed to be $\cap$-flat. The second statement is straighforward.

An easy induction reduces the proof of (c) to the case $s=2$.
Now $A_1,A_2\subseteq A$ so $A_1,A_2\subseteq \twiddle{A}{e}^{[p^e]}$,
and the minimality of $\twiddle{A_1}{e}$ and $\twiddle{A_2}{e}$ now implies
$\twiddle{A_1}{e}, \twiddle{A_2}{e} \subseteq \twiddle{A}{e}$ hence
$\twiddle{A_1}{e}+ \twiddle{A_2}{e} \subseteq \twiddle{A}{e}$.
On the other hand
$$A=A_1+A_2\subseteq \twiddle{A_1}{e}^{[p^e]}+ \twiddle{A_2}{e}^{[p^e]}=
\left(\twiddle{A_1}{e} + \twiddle{A_2}{e}\right)^{[p^e]}$$
and the minimality of
$\twiddle{A}{e}$ implies $\twiddle{A}{e}\subseteq \twiddle{A_1}{e} + \twiddle{A_2}{e}$.
\end{proof}

Notice that if $T$ is a polynomial ring $\mathbb{K}[x_1, \dots, x_n]$ for some field $\mathbb{K}$ of characteristic $p>0$
or a localization of it,
then $T^{1/p^e}$ is a free $T$-module and hence $\cap$-flat.
When $T$ is a power series ring $\mathbb{K}[\![x_1, \dots, x_n]\!]$,
$T^{1/p^e}$ is free $T$-module when $\mathbb{K}^{1/p}$ is finite extension of $\mathbb{K}$, i.e., when $\mathbb{K}$ is \emph {$F$-finite}, but not in general.
However, if the coefficients of a set of generators of the ideal $A\subseteq T$ lie in an $F$-finite field,
the calculation of $\twiddle{A}{e}$ can be carried out over that field.
In the general case we have the following.

\begin{prop}
Let $T=\mathbb{K}[\![x_1, \dots, x_n]\!]$. The $T$-modules $T^{1/p^e}$ are $\cap$-flat for all $e\geq 1$.
\end{prop}
\begin{proof}
It is enough to prove the statement for $e=1$ and we henceforth assume this case.

The fact that $T^{1/p}$ is $T$-flat follows from Theorem 2.1 in \cite{Kunz}.

The rest of this proof
follows the idea described in page 41 of \cite{Hochster-Huneke-1}: we show that if
$\phi: (A,\mathfrak{a}) \rightarrow (B,\mathfrak{b})$ is a flat local map of complete local rings then $B$ is $\cap$-flat over $A$.
Let $N$ be a finitely generated $A$ module and let $\left\{ N_\lambda \right\}_{\lambda\in \Lambda}$ be a set of submodules of $N$.
We show that
$$B \otimes_A \cap_{\lambda\in\Lambda} N_\lambda= \cap_{\lambda\in\Lambda} \left(B \otimes_A N_\lambda\right) .$$
By replacing $N$ with $N/\cap_{\mu\in\Lambda} N_\mu$ and each $N_\lambda$ with its image in $N/\cap_{\mu\in\Lambda} N_\mu$
while using the fact that $B$ is $A$-flat
we may assume that $\cap_{\lambda\in\Lambda} N_\lambda=0$;
after this reduction we need to show that $\cap_{\lambda\in\Lambda} \left(B \otimes_A N_\lambda\right)=0$.

If $\Lambda$ is finite, the result follows directly from the flatness of $B$ so we assume that $\Lambda$ is infinite.
We now reduce to the case where $\Lambda$ is countable by constructing a sequence
$\left\{N_i\right\}_{i\in \mathbb{N}} \subseteq \left\{ N_\lambda \right\}_{\lambda\in \Lambda}$ for which
$\cap_{i\in \mathbb{N}} N_i=0$.
We construct this sequence inductively so that for each $j\geq 0$ there exists an $i_j\geq 1$ such that
$N_1 \cap  \dots \cap N_{i_j} + \mathfrak{a}^j N = \cap_{\lambda\in\Lambda} N_\lambda + \mathfrak{a}^j N $; it is easy to do this for $j=0$ and,
if for some $j\geq 0$
we already defined $N_1 \cap  \dots \cap N_{i_j}$,
we use the fact that the module $N/\mathfrak{a}^{j+1} N$ satisfies DCC to pick
a finite set $N^1, \dots , N^s\subseteq \left\{N_\lambda\right\}_{\lambda\in\Lambda}$ such that
$N_1 \cap  \dots \cap N_{i_j} \cap N^1 \cap  \dots  \cap N^s + \mathfrak{a}^{j+1} N = \cap_{\lambda\in\Lambda} N_\lambda + \mathfrak{a}^{j+1} N$;
we now extend the sequence to $N_1, \dots, N_{i_j}, N^1, \dots, N^s$ and set
$i_{j+1}=i_j+s$.
For the sequence thus constructed we have
$$\cap_{i\in \mathbb{N}} N_i + \mathfrak{a}^{j} N =  \cap_{\lambda\in\Lambda} N_\lambda + \mathfrak{a}^{j} N$$
for all $j\geq 0$ and hence
$\cap_{i\in \mathbb{N}} N_i =  \cap_{\lambda\in\Lambda} N_\lambda =0$.
Assume henceforth that $\Lambda=\mathbb{N}$;
we may replace each $N_i$ with $N_1 \cap \dots \cap N_i$ and assume further that $\left\{N_i\right\}_{i\in \mathbb{N}}$ is decreasing.

We now use Chevalley's Theorem (see Theorem 1 in Chapter 5 of \cite{Northcott}) to deduce that
for all $j>0$ there exists an $i_j$ such that $N_{i_j}\subseteq \mathfrak{a}^j N$.
For all $j\geq 1$ we have
$$B \otimes_A \cap_{i\in \mathbb{N}} N_i \subseteq
B \otimes_A  \mathfrak{a}^j N \subseteq
\mathfrak{b}^j B \otimes_A  N =
\mathfrak{b}^j (B \otimes_A  N)$$
so
$$B \otimes_A \cap_{i\in \mathbb{N}} N_i \subseteq \cap_{j\in \mathbb{N}} \mathfrak{b}^j (B \otimes_A  N)=0 .$$

\end{proof}

Throughout the remainder of this section we will assume that
$T=\mathbb{K}[\![x_1, \dots, x_n]\!]$ or that $T=\mathbb{K}[x_1, \dots, x_n]$ for some field $\mathbb{K}$
of prime characteristic $p$.
We also fix an $e\geq 1$.

Proposition \ref{Proposition: properties of widetilde} reduces the calculation of $\twiddle{A}{e}$ to the case where $A$ is principal,
and this is the content of the next proposition. This proposition has been proved in \cite{ABL} and we reproduce the proof for the reader's convenience.

\begin{prop}\label{Proposition: widetilde of a principal ideal}
Assume that $T$ is free over $T^p$ and let $\mathcal{B}$ be a free basis.
Let $g\in T$ and write $g=\sum_{b\in \mathcal{B}} g_b^{p^e} b$
where $g_b\in T$ for all $b\in \mathcal{B}$.
Then $\twiddle{g T}{e}$ is the ideal generated by
$\left\{g_b \,|\, b\in \mathcal{B} \right\}$.
\end{prop}

\begin{proof}
If $L\subseteq T$ is such that $g=\sum_{b\in \mathcal{B}} g_b^{p^e} b\in L^{[p^e]}$ then we can find
$\ell_1,\dots, \ell_s\in L$ and $r_1, \dots, r_s\in T$ such that
$\sum_{b\in \mathcal{B}} g_b^{p^e} b =\sum_{i=1}^s r_i \ell_i^{p^e}$.
For all $1\leq i\leq s$ we can now write $r_i =\sum_{b\in \mathcal{B}} r_{b,i}^{p^e} b$
where $r_{b,i}\in T$ for all $b\in \mathcal{B}$ and we obtain
$$\sum_{b\in \mathcal{B}} g_b^{p^e} b= \sum_{b\in \mathcal{B}} \left(\sum_{i=1}^s r_{b,i}^{p^e} \ell_i^{p^e}\right) b .$$
Since these are direct sums, we may compare coefficients and deduce that
for all $b\in \mathcal{B}$, $g_b^{p^e}=\sum_{i=1}^s r_{b,i}^{p^e} \ell_i^{p^e}$ hence
$g_b=\sum_{i=1}^s r_{b,i} \ell_i$
and $g_b\in L$. On the other hand, if  $g_b\in L$ for all $b\in \mathcal{B}$ we clearly have
$g=\sum_{b\in \mathcal{B}} g_b^{p^e} b\in L^{[p^e]}$, so we have shown that
$\twiddle{f T}{e}$ is the ideal generated by $\left\{g_b \,|\, b\in \mathcal{B} \right\}$.
\end{proof}

The proposition above translates easily into an algorithm.
Define
$$\Lambda=\left\{ \left(\alpha_1, \dots, \alpha_n\right) \in \mathbb{N}^n \,|\,
 0\leq \alpha_1, \dots, \alpha_n <p^e \right\}$$
and for each $\lambda=\left(\alpha_1, \dots, \alpha_n\right) \in \Lambda$ let $\mathbf{x}^\lambda$ denote
the monomial $x_1^{\alpha_1} \cdot \ldots \cdot x_n^{\alpha_n}$.
Observe next that, if $\Theta$ is a finite basis of $\mathbb{K}$ as a $\mathbb{K}^{p^e}$-vector-space,
$$\mathcal{B}=\left\{ \theta  \mathbf{x}^\lambda \,|\, \theta\in \Theta,\ \lambda\in \Lambda \right\} $$
is a free basis for the $T^{p^e}$-module $T$.

\bigskip
We can now restate one of the statements of Theorem \ref{Theorem: description of nilpotent elements} as follows:
The index of nilpotency of $E_S$, if not zero, is the index at which the descending sequence of ideals
$\left\{ \twiddle{u^{\nu_e}R + I}{e}+I \right\}_{e\geq 1}$ stabilizes.
We can exploit this observation to give a very simple proof, pointed out to me by Gennady Lyubeznik,
of a crucial ingredient used in \cite{ABL}.
Given any $g\in R$, consider the $R[T; f]$-module structure on $E_R$ given by
$\Psi\left( R \xrightarrow[]{g^{p-1}} R\right)$ (here we are taking $I=0$ and $S=R$).
The observation above now implies that the descending chain
$$\left\{ \twiddle{g^{\nu_e(p-1)}R }{e} \right\}_{e\geq 1}= \left\{ \twiddle{g^{(p^e-1)}R }{e} \right\}_{e\geq 1}$$
stabilizes.
If we combine this with Corollary 3.6 in \cite{ABL} we obtain an alternative proof of the fact that
for a power series ring $R$ of prime characteristic, for all nonzero $f\in R$, $1/f$
generates $R_f$ as a $D_R$-module.

More generally, if $G$ is any $m\times m$ matrix with entries in $R$, we may endow
$E_R^m$ with an $R[T; f]$-module structure given by
$\Psi\left( R^m \xrightarrow[]{G} R^m\right)$.
Denote the $(i,j)$ entry of $G^{\nu_e}$ with $g^{(e)}_{i j}$.
It is not hard to see now that
for all $1\leq i, j \leq m$
$$\left\{ \twiddle{g^{(e)}_{i j} R }{e} \right\}_{e\geq 1}$$
is a  descending chain of ideals which stabilizes.

\bigskip
Before proceeding we notice that when $T$ is a polynomial ring and $W\subset T$ is a multiplicative set, Proposition
\ref{Proposition: widetilde of a principal ideal} implies that for any ideal $A\subset T$,
$\twiddle{W^{-1}A}{e}=W^{-1}\twiddle{A}{e}$. Similarly, if $\mathfrak{m}$ is the irrelevant ideal of
$T$ and $\widehat{T}$ denotes the completion of $T$ with respect to  $\mathfrak{m}$, then
$\twiddle{A \widehat{T}}{e}=\twiddle{A}{e} \widehat{T}$.

\bigskip
\begin{defi}\label{definition: star}
Fix  any $u\in T$.
For any ideal $A\subseteq T$  we define a sequence of ideals as follows:
$A_0=A$ and $A_{i+1}=\twiddle{u A_{i}}{e}+{A_i} $ for all $i\geq 0$.
Clearly this sequence is an ascending chain and as $T$ is Noetherian it stabilizes
to some ideal which we denote with $A^{\star^e u}$.
\end{defi}

\begin{prop}\label{Proposition: another property of star}
Fix  any $u\in T$ and let $A\subseteq T$ be an ideal.
If $B\subseteq T$ is an ideal containing $A$ and if $uB \subseteq B^{[p^e]}$ then $A^{\star^e u} \subseteq B$.
\end{prop}
\begin{proof}
Let $\left\{ A_i \right\}_{i=0}^\infty$ be the sequence of ideals as in
Definition \ref{definition: star}. We show by induction that $A_i\subseteq B$ for all $i\geq 0$.
Since $A_0=A$ and $A\subseteq B$ the claim is true for $i=0$; assume that $i\geq 0$ and that $A_i\subseteq B$.
Now $u A_i \subseteq u B \subseteq B^{[p^e]}$ and the minimality of
$\twiddle{u A_{i}}{e}+{A_i}$ now implies that
$A_{i+1}=\twiddle{u A_{i}}{e}+{A_i} \subseteq B$.
\end{proof}

The regular local ring $R$ at the focus of this paper is a power series hence
it is of the form  $\mathbb{K}[\![x_1, \dots, x_n]\!]$
for some field $\mathbb{K}$ of prime characteristic $p$.
When $A$ is expanded from the polynomial ring $T=\mathbb{K}[x_1, \dots, x_n]$ and $u\in T$
we want to compute $A^{\star^e u}$ by performing calculations in $T$ rather than $R$.
Proposition  \ref{Proposition: polynomial star closure} below shows how to do that.

\begin{prop}\label{Proposition: polynomial star closure}
Let $A$ be an ideal of $T=\mathbb{K}[x_1, \dots, x_n]$ and let $u\in T$.
We have $(A R)^{\star^e u}=\left(A^{\star^e u}\right) R$.
\end{prop}
\begin{proof}
Let $\{B_i\}_{i\geq 0}$ and $\{C_i\}_{i\geq 0}$ be the sequences introduced in Definition \ref{definition: star}
whose stable values are $(A R)^{\star^e u}$ and $A^{\star^e u}$, respectively. We will show
that $B_i=C_i R$ for all $i\geq 0$ using induction on $i$.

First, $C_0 R =A R= B_0$, so assume that $i>0$ and that $B_{i-1}=C_{i-1} R$.
Now notice that since
$u C_{i-1} \subseteq C_i^{[p^e]}$ and $C_{i-1}\subseteq C_i$
we have
$u B_{i-1}=u C_{i-1}R  \subseteq C_i^{[p^e]}R = \left(C_i R\right)^{[p^e]}$ and $B_{i-1}=C_{i-1}R\subseteq C_i R$
so the minimality of $B_i$ implies that $B_i\subseteq C_i R$.
On the other hand,
$u C_{i-1} R=u B_{i-1} \subseteq B_i^{[p^e]}$ implies that
$u C_{i-1} =u C_{i-1} R \cap T \subseteq B_i^{[p^e]} \cap T=\left(B_i \cap T\right)^{[p^e]}$
and
$C_{i-1} R=B_{i-1} \subseteq B_i$ implies that
$C_{i-1}= C_{i-1} R\cap T\subseteq B_i\cap T$
and the minimality of $C_i$ implies that
$C_i\subseteq B_i \cap T$ hence
$C_i R \subseteq (B_i \cap T)R\subseteq B_i$.
\end{proof}

\section{$E_S$-ideals and special $E_S$-ideals}

Following \cite{Sharp}, we call an ideal $K\subseteq S[T;f]$ a \emph{graded two-sided ideal} if
$\displaystyle K=\bigoplus_{i=0}^\infty  K_i T^i$ for ideals
$K_0, K_1, \dots$ of $S$. An important example is $K=L S[T; f]$ for some ideal $L\subseteq S$.
Let $G$ be an $S[T;f]$-module.
An $S[T;f]$-submodule $M\subseteq G$ is a \emph{special annihilator submodule}
if $M=\Ann_G K$ for some graded two-sided ideal $K\subseteq S[T;f]$.

For any $S[T;f]$-submodule $M\subseteq G$ we define the \emph{graded annihilator of $M$,} denoted
$\grAnn_{S[T;f]} M$,
to be the largest graded two-sided ideal contained in
$\Ann_{S[T;f]} M$.

We call an ideal $L\subseteq S$ a \emph{$G$-special ideal} whenever $L S[T;f]$ is the graded annihilator of some
$S[T;f]$-submodule $M\subseteq G$, in which case $L S[T;f]=\grAnn \left(\Ann_G L S[T;f]\right)$
(cf. Lemma 1.7 in \cite{Sharp};
notice that we extended slightly the definition of special ideals to the case where $G$ is not necessarily
$T$-torsion-free).

\begin{prop}\label{E-special is E-ideal}
Assume that $R$ is complete and that $E_S$ is $T$-torsion free.
An ideal $L\subseteq R$ which contains $I$ is an $E_S$-ideal if and only if $L S$ is $E$-special.
\end{prop}

\begin{proof}
Assume first that $L$ is $E_S$-special, i.e.,
$L S[T;f]=\grAnn N$ for some $S[T;f]$-submodule $N$ of $E_S$, and since $E_S$ is assumed to be $T$-torsion-free,
we have $\grAnn N= (0 :_R N) S[T;f]$ (cf.~Definition 1.10 in \cite{Sharp}).
We can also write $N=\Ann_{E_S} L^\prime$ for some $E_S$-ideal $L^\prime$ and
\begin{eqnarray*}
L S[T;f] & = & \left(0 :_R \Ann_{E_S} L^\prime\right) S[T;f]\\
&=& \left(0 :_R  \left(R/L^\prime\right)^\vee\right) S[T;f]\\
&=& \left(0 :_R  R/L^\prime\right) S[T;f]\\
&=& L^\prime  S[T;f]
\end{eqnarray*}
so $L=L^\prime$ and is an $E_S$-ideal.

If, on the other hand,
$L$ is an $E_S$-ideal, i.e., if $\Ann_{E_S} L = \Ann_{E_S} L S[T; f]$, then
\begin{eqnarray*}
\grAnn \Ann_{E_S} L S[T;f] & = & ( 0 :_R \Ann_{E_S} L S[T;f] ) S[T;f]\\
& = & ( 0 :_R \Ann_{E_S} LS  ) S[T;f]\\
& = & \left( 0 :_R  \left(R/ L\right)^\vee  \right) S[T;f]\\
& = & ( 0 :_R R/L  ) S[T;f] \\
& = & L S[T;f]
\end{eqnarray*}
and so $L$ is $E_S$-special.

\end{proof}

\section{The $S[T;f]$-module structure of $\H^{\dim S}_{\mathfrak{m} S}(S)$ and the induced structure on $E_S$}
\label{The S[T;f]-module structure of H}

In what follows we describe a natural $S[T;f]$-module structure on
$\H^{\dim S}_{\mathfrak{m} S}(S)$ and show how this induces an $S[T;f]$-module structure on
$E_S$;
the following section will describe its relevance to test-ideals.

We shall assume henceforth that $S$ is Cohen-Macaulay with canonical module
$\omega\subseteq S$.

The short exact sequence
$0\rightarrow \omega \rightarrow S \rightarrow S/\omega \rightarrow 0$
yields a surjection
$\H^{\dim S}_{\mathfrak{m} S}(\omega) \twoheadrightarrow \H^{\dim S}_{\mathfrak{m} S}(S)$;
for any system of parameters $x_1, \dots, x_d$ for $S$ this
map can also be described as the map
$$\underrightarrow{\lim}_{i\geq 0} \frac{\omega}{x_1^i \omega + \dots + x_d^i \omega} \rightarrow
\underrightarrow{\lim}_{i\geq 0} \frac{S}{x_1^i S + \dots + x_d^i S} $$
induced by the maps
$$\frac{\omega}{x_1^i \omega + \dots + x_d^i \omega} \rightarrow
\frac{S}{x_1^i S + \dots + x_d^i S}$$
given by $a+ \left(x_1^i \omega + \dots + x_d^i\omega \right)  \mapsto a+ \left(x_1^i S + \dots + x_d^i S\right)$.
The natural action of Frobenius on $\H^{\dim S}_{\mathfrak{m} S}(S)$
given by
$$f\left( a+ \left(x_1^i S + \dots + x_d^i S\right) \right)=
a^p+ \left(x_1^{i p} S + \dots + x_d^{i p} S\right) $$
now lifts to an action on
$\displaystyle E_S \cong \underrightarrow{\lim}_{i\geq 0} \frac{\omega}{x_1^i \omega + \dots + x_d^i \omega}$
given by
$$f\left( a+ \left(x_1^i \omega + \dots + x_d^i \omega\right)  \right)=
a^p+ \left(x_1^{i p} + \dots + x_d^{i p} \omega\right)$$
and this $S[T;f]$ module structure on $E_S$ clearly makes the surjection
$E_S \twoheadrightarrow \H^{\dim S}_{\mathfrak{m} S}(S)$ described above into
a $S[T;f]$-linear map.

If we apply $\Delta$ to the $S[T;f]$-linear
surjection $E_S \twoheadrightarrow \H^{\dim S}_{\mathfrak{m} S}(S)$
and  identify $E_S^\vee$ with $R/I$ we obtain the following commutative diagram with exact rows
\begin{equation}\label{CD2}
\xymatrix{
0 \ar@{>}[r]^{} &
\displaystyle\H^{\dim S}_{\mathfrak{m} S}(S)^\vee \ar@{>}[r]^{} \ar@{>}[d]^{u} &
\displaystyle\frac{R}{I} \ar@{>}[d]^{u} \\
0 \ar@{>}[r]^{} &
\displaystyle \displaystyle F_R\left(\H^{\dim S}_{\mathfrak{m} S}(S)^\vee\right) \ar@{>}[r]^{} &
\displaystyle \displaystyle\frac{R}{I^{[p]}}
}
\end{equation}
where $u\in R$, the second vertical map is multiplication by $u$ and the first vertical map
is given by restriction of the second, i.e., also by multiplication by $u$.
We deduce that under the identification of $E_S^\vee$ with $R/I$,
$\displaystyle\H^{\dim S}_{\mathfrak{m} S}(S)^\vee$ is identified with
$J/I$ for some ideal $ J\subseteq R$ containing $I$.
This ideal $J$ must then satisfy $u J \subseteq J^{[p]}$.

Our next step is to compute $J$ and $u$ effectively.
Let $\Omega$ be the full pre-image of $\omega$ in $R$.
Working over $R$, the surjection
$E_S=\H^{\dim S}_{\mathfrak{m} S}(\omega) \twoheadrightarrow \H^{\dim S}_{\mathfrak{m} S}(S)$
can be written as
$E_S=\H^{\dim S}_{\mathfrak{m}} \left(\Omega/I \right) \twoheadrightarrow
\H^{\dim S}_{\mathfrak{m}}(R/I)$.
Write $\delta=\dim R-\dim S$;
recall that Local Duality states that the functors $\H^{\dim S}_{\mathfrak{m}} \left( - \right)$
and
$\Ext_R^\delta\left( -, R \right)^\vee$
are isomorphic and the surjection above is induced by applying either of these functors
to the inclusion $\omega\subseteq S$.
Applying the latter and a further application of
$(-)^\vee$ yields the injection
$\Ext_R^\delta\left( R/I, R \right) \subseteq \Ext_R^\delta\left( \Omega/I, R \right)$
and so
$J/I\cong \Ext_R^\delta\left( R/I, R \right)$.
This $\Ext$-module can  be computed effectively and $J$ can be recovered by computing
a minimal presentation of this module.


To find the map $u$ in (\ref{CD2}) we use the fact that
$\mathcal{F} \left(\H^{\dim S}_{\mathfrak{m} S}(S)\right)$
is the $R$-algebra with one generator corresponding to the $S[T;f]$-module structure defined above (cf.~Example 3.7 in \cite{Lyubeznik-Smith}).
Hence the $S$-linear maps
$\Ext_R^\delta\left( R/I, R \right) \rightarrow F_R \left(\Ext_R^\delta\left( R/I, R \right) \right)$
form a rank-one $S$-free module and the generator $u$ of this free module
can be computed explicitly from the generator of
$$\frac{\left(I^{[p]} :_R I\right) \cap \left(J^{[p]} :_R J\right)}{I^{[p]}} $$
(cf.~Chapter 3 of \cite{Blickle} and section \ref{Section: Frobenius maps on injective hulls}).

\section{The computation of parameter-test-ideals}
Throughout this section we will assume that $S=R/I$ is Cohen-Macaulay with canonical module $\omega\subseteq S$.
We shall write $H=\H^{\dim S}_{\mathfrak{m} S}(S)$ and we will assume that $E_S$ is $T$-torsion-free.
This last assumption implies that
$H$, being a quotient of $E_S$ by a special annihilator submodule, is also $T$-torsion-free (cf.~Lemma 3.1 in \cite{Sharp}).
Recall now that Corollary 4.6 in \cite{Sharp} now states that the parameter-test-ideal of $S$
is the smallest $H$-special ideal of $S$ of positive height.
In this section we relate the $H$-special ideals to $E_S$-special ideals
and describe an algorithm for computing the parameter-test-ideal of $S$.

First we note the following.

\begin{thm}\label{Theorem: spcial ideals in quotient}
Assume that $E_S$ is $T$-torsion-free and write
$H:=\displaystyle \H^{\dim S}_{\mathfrak{m} S}(S)\cong \frac{E_S}{\Ann_{E_S} J}$
where $J\subseteq R$ is an $E_S$-ideal.
The $H$-special ideals are
$$ \left\{ (L : J) \,|\, L\subseteq R \text{ is a } E_S \text{-special ideal contained in } J \right\} .$$
\end{thm}
\begin{proof}
This follows from Proposition \ref{E-special is E-ideal} and  Proposition 3.3 in \cite{Sharp}.
\end{proof}

We are now ready for the main theorem in this section.
\begin{thm}
Assume that $E_S$ is $T$-torsion-free.
Let $c\in R$ be such that its image in $S$ is a parameter test element.
The parameter test ideal $\overline{\tau}$ of $S$ is given by
$ ((cJ+I)^{\star u} :_R J)S$.
\end{thm}
\begin{proof}
Notice that $(cJ+I)^{\star u}$ is an $E_S$-ideal and that since $c\in ((cJ+I)^{\star u} : J)$, we have
\goodbreak $\height ((cJ+I)^{\star u}: J)S >0$.
Now
\begin{eqnarray*}
\overline{\tau} & = & \cap \{K \,|\, K\subset S \text{ is a $H$-special ideal}, \height K > 0 \}\\
&=& \cap \{(L:_R J) \,|\, L\subset J \text{ is an $E_S$-ideal}, \height (L:_R J)S > 0 \}\\
&=& \big(\cap \{L \,|\, L\subset J \text{ is an $E_S$-ideal}, \height (L:J) S > 0 \} : J\big)
\end{eqnarray*}
so we see that $\overline{\tau} \subseteq ((cJ+I)^{\star u} : J)$.

Also, $c\in \overline{\tau}$ hence
$cJ\subseteq L$ for all  $E_S$-ideals $L$ for which $\height (L:J)S > 0$ and
Proposition \ref{Proposition: another property of star}
implies that $(cJ+I)^{\star u} \subseteq L$ and hence that
$\left( (cJ+I)^{\star u} : J \right)  \subseteq (L: J)$
for all  $E_S$-ideals $L$ for which  $\height (L:J) > 0$.
We conclude that $\left( (cJ+I)^{\star u} : J \right)\subseteq \overline{\tau}$.
\end{proof}

In the case where $E_S$ is $T$-torsion-free, if we are given \emph{one} parameter-test-element,
we can now \emph{compute} the \emph{entire} parameter test ideal of $S$
as follows.

\begin{enumerate}
    \item Find the element $u\in R$ as described in section \ref{The S[T;f]-module structure of H}
    and use Theorem \ref{Theorem: description of nilpotent elements} to determine whether
    $E_S$ is $T$-torsion-free. If $E_S$ is $T$-torsion-free proceed as follows.
    \item Find the ideal $I\subseteq J\subseteq R$  as described
    in section \ref{The S[T;f]-module structure of H}.
    \item Given one parameter test element $c$, compute $L=(cJ+I)^{\star u}$ as
    described in section \ref{Section: The star-closure}.
    \item The parameter test ideal of $S$ is $(L :_R J)S$.
\end{enumerate}

We also note that the verification of whether $E_S$ is $T$-torsion-free is also algorithmic:
the proof of Theorem \ref{Theorem: description of nilpotent elements} shows that
$E_S$ is $T$-torsion-free if and only if $I_1(u)+I=R$.

\section{Applications and examples}

A particularly simple instance of the results of the previous chapters is the case where
$S$ is a complete intersection, i.e., the case where $I$ is generated by a regular sequence $u_1, \dots, u_s\in R$.
Now $S=R/I$ is Gorenstein, $E_S=\HH^{\dim S}_{\mathfrak{m} S}(S)$ (so the surjection described in the beginning of
section \ref{The S[T;f]-module structure of H} is an equality) and
$\Delta(E_S)=\left( R/I \rightarrow R/I^{[p]} \right)$ is given by multiplication by $u=(u_1 \cdot \ldots \cdot u_s)^{p-1}$
whose image in the $S$-module $(I^{[p]} :_R I)/I^{[p]}$ generates it.

We call $S$
\emph{$F$-injective} if the natural Frobenius map on $E_S=\HH^{\dim(S)}_{\mathfrak{m}S}(S)$ is injective, i.e.,
if $\Nil(E_S)=0$. We now recover Fedder's Criterion (Proposition 2.1 in \cite{Fedder})
which states that, with $S$ as in the previous paragraph,
$S$ is $F$-injective
if and only if $u\notin \mathfrak{m}^{[p]}$ where $u=(u_1 \cdot \ldots \cdot u_s)^{p-1}$.
The crucial fact here is that $\Delta(E_S)$ is the map
$R/I \xrightarrow[]{u} R/I^{[p]}$.
As in the proof of Theorem \ref{Theorem: description of nilpotent elements}
consider $N_1=\{ m\in E_S \,|\, T m=0\}$ and write $\Delta(N_1)=(R/L \xrightarrow[]{u} R/L^{[p]})$
for some $E_S$-ideal $L$.
We saw that this map is the zero map, i.e., $u\in L^{[p]}$.
Fedder's condition $u\notin \mathfrak{m}^{[p]}$ is equivalent to the non-existence of a proper ideal $L\subset R$
for which $u\in L^{[p]}$ so it implies that $N_1=\Ann_{E_S} R= 0$.
If, on the other hand,
$u\in \mathfrak{m}^{[p]}$ then $u  \mathfrak{m} \subseteq  \mathfrak{m}^{[p]}$, so $\mathfrak{m}$
is an $E_S$ ideal and since the map
$R/\mathfrak{m} \xrightarrow[]{u} R/\mathfrak{m}^{[p]}$ is the zero map,
$T \Psi\left(R/\mathfrak{m} \xrightarrow[]{u} R/\mathfrak{m}^{[p]}\right)=0$
and $S$ is not $F$-injective.

\bigskip
We now describe a specific calculation performed using the methods in the previous sections.
All calculations described below were performed with Macaulay2 \cite{M2}.

Let $\mathbb{K}$ be the field of two elements, $R=\mathbb{K}[x_1, x_2, x_3, x_4, x_5]$,
let $I$ be the ideal of $R$ generated by the $2\times 2$ minors of
$$\left( \begin{array}{llll} x_1 & x_2 & x_2  & x_5 \\ x_4 & x_4 & x_3 & x_1 \end{array} \right)$$
and let $S=R/I$. This quotient is reduced,  $2$-dimensional, Cohen-Macaulay and of Cohen-Macaulay type 3;
we produce a canonical module by computing
$$\Ext^3_R(S,R)\cong \coker
\left( \begin{array}{llllllll}
{x}_{2}&       {x}_{1}&       0&       0&       {x}_{3}+{x}_{4}&       {x}_{4}&       {x}_{5}&       {x}_{4}\\
       0&       0&       {x}_{3}&       {x}_{4}&       0&       0&       {x}_{1}&       0\\
       {x}_{5}&       {x}_{5}&       {x}_{5}&       {x}_{5}&       0&       {x}_{2}&       0&       {x}_{1}
\end{array} \right);
$$
this is isomorphic to the ideal $\omega\subset S$ which is the image in $S$ of the ideal $\Omega\subset R$
generated by $x_1, x_4, x_5$.

We now take $J=\Omega$ and compute the generator $u$ of the $S$-module
$$\frac{\left(I^{[2]} :_R I\right) \cap \left(J^{[2]} :_R J\right)}{I^{[2]}} ;$$
this turns out to be
$$u={x}_{1}^{3} {x}_{2} {x}_{3}+{x}_{1}^{3} {x}_{2} {x}_{4}+{x}_{1}^{2} {x}_{3} {x}_{4} {x}_{5}+{x}_{1} {x}_{2} {x}_{3} {x}_{4} {x}_{5}+{x}_{1} {x}_{2} {x}_{4}^{2}
       {x}_{5}+{x}_{2}^{2} {x}_{4}^{2} {x}_{5}+{x}_{3} {x}_{4}^{2} {x}_{5}^{2}+{x}_{4}^{3} {x}_{5}^{2} .$$
We compute $\twiddle{u^{\nu_1}R}{1}=R$ hence $E_S$ is $T$-torsion free.
Now the parameter-test-ideal $\tau$ is computed as $\left( (cJ+I)^{\star u} : J \right)$
where $c$ is randomly chosen to be in the defining ideal of the singular locus of $S$ and not in a minimal prime of $I$. This calculation yields
$\tau=(x_1,x_2, x_3+x_4, x_4x_5)R$ and we deduce that $S$ is not $F$-rational.


\begin{thebibliography}{HKSY}

\bibitem[ABL]{ABL}
J.~Alvarez-Montaner, M.~Blickle, G.~Lyubeznik.
\emph{Generators of D-modules in positive characteristic.}
Mathematical Research Letters,  {\bf 12}  (2005),  no.~4, pp.~459--473.

\bibitem[BK]{Brenner-Katzman}
H.~Brenner and M.~Katzman.
\emph{On the arithmetic of tight closure.}
Journal of the AMS {\bf 19}  (2006),  no.~3, pp.~659--672.

\bibitem[B]{Blickle}
M.~Blickle.
\emph{The intersection homology D--module in finite characteristic.}
PhD thesis, University of Michigan, 2001.
ArXiV math.AG/0110244.



\bibitem[F]{Fedder}
R.~Fedder.
\emph{$F$-purity and rational singularity.}
Transactions of the AMS, {\bf 278} (1983), no.~2, pp.~461--480.

\bibitem[GS]{M2}
D.~R.~Grayson and M.~E.~Stillman.
\emph{Macaulay 2, a software system for research in algebraic geometry},
available at http://www.math.uiuc.edu/Macaulay2/.




\bibitem[H]{Huneke}
C.~Huneke.
\emph{Tight closure and its applications.}
CBMS Regional Conference Series in Mathematics, 88.
American Mathematical Society, Providence, RI, 1996.

\bibitem[HS]{Hartshorne-Speiser} R.~Hartshorne and R.~Speiser.
\emph{Local cohomological dimension in characteristic $p$,} Ann.~of Math.~
\textbf{105} (1977), pp.~45--79.

\bibitem[HKSY]{Huneke-Katzman-Sharp-Yao} C.~Huneke, M.~Katzman, R.~Y.~Sharp and Y.~Yao.
\emph{Frobenius test exponents for parameters ideals in generalized Cohen-Macaulay local rings,}
Journal of Algebra,  {\bf 305} (2006),  pp.~516-539.



\bibitem[HH1]{Hochster-Huneke-0}
M.~Hochster and C.~Huneke.
\emph{Tight closure, invariant theory, and the Brian\c{c}on-Skoda theorem.}
Journal of the AMS {\bf 3} (1990), no. 1, 31--116.


\bibitem[HH2]{Hochster-Huneke-1}
M.~Hochster and C.~Huneke.
\emph{$F$-regularity, test elements, and smooth base change.}
Transactions of the AMS, {\bf 346} (1994), no.~1, pp.~1--62.



\bibitem[K]{Kunz}
E.~Kunz.
\emph{Characterizations of regular local rings for characteristic $p$.}
American Journal of Mathematics {\bf 91} (1969), pp.~772--784.

\bibitem[KS] {Katzman-Sharp} M.~Katzman and R.~Y.~Sharp.
\emph{Uniform behaviour of the Frobenius closures of ideals generated by regular sequences,}
Journal of  Algebra,  \textbf{295} (2006) 231--246.

\bibitem[L]{Lyubeznik}
G.~Lyubeznik.
\emph{$F$-modules: applications to local cohomology and $D$-modules in characteristic $p>0$.}
J.~Reine Angew.~Math.~{\bf 491} (1997), pp.~65--130.

\bibitem[LS]{Lyubeznik-Smith}
G.~Lyubeznik and K.~E.~Smith.
\emph{On the commutation of the test ideal with localization and completion.}
Transactions of the AMS  {\bf 353} (2001), no.~8, pp.~3149--3180.

\bibitem[N]{Northcott}
D.~G.~Northcott.
\emph{Ideal theory.}
Cambridge Tracts in Mathematics and Mathematical Physics, {\bf 42}. Cambridge University Press, 1953.

\bibitem[S]{Sharp}
R.~Y.~Sharp.
\emph{Graded annihilators of modules over the Frobenius skew polynomial ring, and tight closure.}
Transactions of the AMS  {\bf 359}  (2007),  no.~9, pp.~4237--4258

\bibitem[SV]{Sharpe-Vamos}
D.~W.~Sharpe and P.~V{\'a}mos.
\emph{Injective modules.}
Cambridge Tracts in Mathematics and Mathematical Physics {\bfseries 62}.
Cambridge University Press, London-New York, 1972.




\bibitem[Sm]{Smith2}
K.~E.~Smith.
\emph{Test ideals in local rings.}
Transactions of the AMS {\bf 347} (1995), no.~9, pp.~3453--3472.



\end{thebibliography}
\end{document}